\documentclass[reqno,11pt]{amsart} 

\usepackage[colorlinks=true,urlcolor=blue,
citecolor=red,linkcolor=blue,linktocpage,pdfpagelabels,
bookmarksnumbered,bookmarksopen]{hyperref}

\usepackage[left=2.6cm,right=2.6cm,top=2.9cm,bottom=2.9cm]{geometry}

\usepackage{amsmath,amssymb,latexsym,soul,cite,mathrsfs}
\numberwithin{equation}{section}  % Equation within section - Works with amsmath
\usepackage{color,enumitem,graphicx}

\usepackage{cancel}

\usepackage[colorlinks=true,urlcolor=blue,
citecolor=red,linkcolor=blue,linktocpage,pdfpagelabels,
bookmarksnumbered,bookmarksopen]{hyperref}
\usepackage[english]{babel}

\newtheorem{theorem}{Theorem}[section]
\newtheorem{proposition}[theorem]{Proposition}
\newtheorem{corollary}[theorem]{Corollary}
\newtheorem{lemma}[theorem]{Lemma}
\DeclareMathAlphabet{\mathpzc}{OT1}{pzc}{m}{it}

\newtheorem{remark}[theorem]{Remark}

\title[Geometry of Einstein-type manifolds with boundary]{Geometry of Einstein-type manifolds with boundary}

\author[M. Andrade]{Maria Andrade}
\address[M. Andrade]{Departamento de Matemática, 
 Universidade Federal de Sergipe, Brazil
	\newline\indent 
}
\email{\href{mailto: maria@mat.ufs.br}{ maria@mat.ufs.br}}

\thanks{The author was partially supported by Brazilian National Council for Scientific and Technological Development (CNPq Grants 403349/2021-4, 408834/2023-4 and 403869/2024-2) and FAPITEC/SE/Brazil 019203.01303/2024-1.}


\begin{document}

\maketitle
\hspace{5.0cm} To the memory of Celso Viana
\begin{abstract}
In this article, we consider Einstein-type manifolds with boundary which generalizes important geometric equations, like static vacuum and static perfect fluid. We investigate some geometric inequalities for those manifolds. Then, we established boundary estimates in terms of the first eigenvalue of the Jacobi operator and another one related to the Brown-York mass.
\end{abstract}

% \tableofcontents

\section{Introduction}\label{intro}
A Riemannian manifold $(M^n,g),$ $n\geq 3$, will be called an Einstein-type manifold, if there exists smooth functions $f,h:M\to\mathbb{R}$ such that
\begin{eqnarray}\label{eq0}
fRic=\nabla^2f+hg,
\end{eqnarray}
where $\nabla^2f$ stands for the Hessian of $f$ and $Ric$ is the Ricci tensor for $M$, $f> 0$ in int(M) and $f=0$ on $\partial M$ (see \cite{leandro2021vanishing}).
We observe that the function $h$ is determined by $f$ and the scalar curvature of $g$, denoted by $R$. In fact, if we trace the equation \eqref{eq0}, then we obtain
\begin{eqnarray}\label{eq1}
fR=\Delta f+nh,
\end{eqnarray}
here $\Delta$ is the Laplacian operator. 

Some well known structural equations are examples of Einstein-type manifolds. For instance, if $f$ is a constant positive, $\partial M=\emptyset,$ then $h=fR/n$ and $(M^n,g)$ is an Einstein manifold. If $h=0$ and $R=0,$ we obtain the static vacuum Einstein equation. If $h=Rf/(n-1)$, we have the vacuum static equation. If we consider $h=Rf/(n-1)+k/(n-1),$ where $k$ is a constant, then we obtain the $V$-static equation. Finally, if $h=(\mu-\rho)f/(n-1),$ where $\mu=R/2$ and $\rho$ are, respectively, the density and the pressure smooth functions, we obtain the static perfect fluid (see \cite{andrade2024some} and references therein).

We recall that there exists a classic conjecture called cosmic no-hair conjecture, proposed by Boucher et al. \cite{boucher1984uniqueness}, which is rewritten as follow:

\textit{The only $n$-dimensional compact static $(M^n, g,f)$ with positive scalar curvature and connected boundary $\partial M$ is given by a round hemisphere $\mathbb{S}_+^n,$ where the function $f$ is taken as the height function.}

Since geometric inequalities are fundamental objects of study in geometry, Boucher-Gibbons-Horowitz \cite{boucher1984uniqueness} and Shen \cite{shen1997note} gave an important answer to this question. More precisely, they proved the following result:
\begin{theorem}[Boucher-Gibbons-Horowitz \cite{boucher1984uniqueness}, Shen \cite{shen1997note}] Let $(M^3,g,f)$ be a compact oriented static triple with connected boundary $\partial M$ and scalar curvature $6.$ Then $\partial M$ is a two-sphere whose area satisfies the inequality
$$|\partial M|\leq 4\pi$$
which equality holding if and only if $M^n$ is isometric to the round hemisphere.
\end{theorem}

Motivated by the last result among others, the named author with the de Melo \cite{andrade2024some} obtained a natural extension for Einstein-type manifold. To be precise, they proved the following result.

\begin{theorem}[\cite{andrade2024some}]\label{yamabe_vol}
Let $(M^n,g,f,h)$ be a compact oriented Einstein-type manifold with oriented connected Einstein boundary $\partial M$  with  $R^{\partial M}>0$ and $R_{min}=n(n-1),$ where $R_{min}$ is the minimum value of $R$ on $M^n$. Suppose that $h\geq fR/n.$ Then,
\begin{eqnarray}\label{yamabe_vole}
|\partial M|\leq \omega_{n-1},
\end{eqnarray}
here $\omega_{n-1}$ denotes the volume of the round unit sphere $\mathbb{S}^{n-1}$. The equality holds in \eqref{yamabe_vole} if and only if $(M^n,g)$ is an Einstein manifold, in this case $(M^n,g)$ is isometric to a standard round hemisphere.
\end{theorem}

To proceed we need to recall some basic facts about the Yamabe constant. We consider $(\Sigma^n,g)$ be a closed smooth manifold of dimension $n$ and we denote by $[g]$ the conformal class of a Riemannian metric $g$ on $\Sigma^n,$ i .e., we say that $\tilde{g}\in [g],$ if there exists a positive smooth function $F$ such that $\tilde{g}=Fg.$ The Yamabe constant of $[g]$ denoted by $Y(M,[g])$ is defined by
 \[
Y(\Sigma,[g])=\inf_{\tilde{g}\in [g]}\frac{\int R_{\tilde{g}}dV_{\tilde{g}}}{\left(\int_{\Sigma}dV_{\tilde{g}}\right)^{\frac{n-2}{n}}}.
 \]

In particular, it is well known that the Yamabe constant of a standard sphere $\mathbb{S}^n$ is given by $Y(\mathbb{S}^n,[g_{can}])=n(n-1)\omega_n^{2/n}.$ Thus, under hypothesis of Theorem \ref{yamabe_vol}, we obtain that 
$$|\partial M|\leq \frac{Y(\mathbb{S}^{n-1},g_{can})^{\frac{n-1}{2}}}{(n-1)(n-2)}$$
and the equality holds if and only if $M$ is isometric to the round hemisphere.

It is possible to prove that the boundary $\partial M$ of an Einstein-type manifold is totally geodesic (Section \ref{prel}). This motivated us to give an estimate for the area of the boundary $\partial M$ in terms of the eigenvalue of Jacobi operator of $\partial M.$ We consider $\varphi\in C^{\partial M},$ then the Jacobi operator, denoted by $J$, is defined by
 $$J(\varphi)=\Delta_{\partial M}\varphi+(Ric(\nu,\nu)+|II|^2)\varphi,$$
 where $\Delta_{\partial M}$ stands for the Laplacian operator on $\partial M,$ $Ric(\nu,\nu)$ is the Ricci curvature of $M$ in the direction of the outward unit normal vector field $\nu$ and $II$ is the second fundamental form of $\partial M.$ The first eigenvalue of the Jacobi operator $J$, is defined by
 \[
 \lambda_1=\inf_{\varphi\neq 0}\frac{-\displaystyle\int_{\partial M} J(\varphi)dS_g}{\displaystyle\int_{\partial M}\varphi^2dS_g}
 \]
 
 and $J(\varphi)=-\lambda_1\varphi.$
Inspired by ideas in (\cite{barros2019rigidity}, Theorem 1.10) and Theorem 3 in \cite{costa2023geometry}, we will prove the following result.

\begin{theorem}\label{estlam} Let $(M^n, g,f,h)$ be a compact Einstein-type manifold with connected Einstein boundary and $R^{\partial M}>0.$ Then,
\begin{eqnarray}\label{lambda_1}
\lambda_1\leq \frac{1}{2}(Y(\mathbb{S}^{n-1}, g_{can})|\partial M|^{\frac{-2}{n-1}}-R_{min}),
\end{eqnarray}
where $R_{min}=\text{min}\{R(p);p\in \partial M\}.$
In addition, if equality holds in \eqref{lambda_1}, then $(\partial M, g_{\partial M})$ is isometric, up to scaling, to the standard sphere $\mathbb{S}^{n-1}.$
\end{theorem}

  Barros and Silva \cite{barros2019rigidity} used the first eigenvalue of the Jacobi operator to proved an estimate for the area of the boundary of a triple static. A similar result was obtained by Costa et al. \cite{costa2023geometry} for a static perfect fluid fluid space-time under certain conditions. Recently, Diógenes et. al \cite{diogenes2024geometric} established a boundary estimate for compact quasi-Einstein manifold with boundary in terms of the first eigenvalue of the Jabobi operator $\lambda_1$. Here, inspired by \cite{barros2019rigidity, costa2023geometry} and ideas outlined in (\cite{diogenes2024geometric}, Theorem 3), we shall establish a similar estimate for the area of the boundary $\partial M$ of a compact Einstein-type manfold. More precisely, we have the following result.

\begin{theorem}\label{boundlambda1}
Let $(M^n,g,f,h)$, $n\geq 3,$ be an Einstein-type manifold with connected Einstein boundary,
$R_{min} >0$ and $R^{\partial M}\geq R_{max}>0.$ Then, one has
\begin{eqnarray}\label{estboun}
|\partial M|\leq \left(\frac{(n-1)(n-2)}{2\lambda_1+(n-1)}\right)^{\frac{n-1}{2}}\omega_{n-1},
\end{eqnarray}
where $\omega_{n-1}$ denotes the volume of the round unit sphere $\mathbb{S}^{n-1}$ and $R_{min},\ R_{max}$ are the minimum and maximum value of $R$ on $\partial M,$ respectively. Moreover, if the equality holds in \eqref{estboun}, then $\partial M$ is isometric to the round sphere $\mathbb{S}^{n-1}$.
\end{theorem}

Before presenting our next result, it is necessary to recall the definition of the Riemannian Brown-York mass. Let $(M^n,g)$ be a Riemannian manifold and consider $\Sigma$ as a connected hypersurface in $(M^n,g)$ such that $(\Sigma, g|_{\Sigma})$ can be embedded in $\mathbb{R}^n$ as a convex hypersurface. The Riemannian Brown-York mass denoted by $m_{BY}$ of $\Sigma$ with respect to $g$ is given by
$$m_{BY}(\Sigma)=\displaystyle\int_{\Sigma}(H_0-H_g)dS_g,$$
here $H_0$ and $H_g$ are the mean curvature of $\Sigma$ as hypersurface of $\mathbb{R}^n$ and $M,$ respectively, and $dS_g$ is the volume element of $\Sigma$ induced by $g.$ Yuan, inspired by Riemannian Penrose inequality, proved a boundary estimate for static spaces in terms of the Riemannian Brown-York mass. In \cite{costa2023geometry} was established a sharp boundary estimate for compact static perfect fluid space-time with (possibly disconneted) boundary in terms of the Riemannian Brown-York mass. Inspired by these results, the following Theorem holds.

\begin{theorem}\label{massc}
Let $(M^n, g,f,h)$ be a compact type-Einstein manifold, $n\geq 3$ with (possibly disconnected) boundary and positive scalar curvature satisfying $h\leq R_gf/(n-1)$. Suppose that each boundary component $(\partial M_i,g)$ can be isometrically embedded in $\mathbb{R}^n$ as a convex hypersurface. Then we have
$$|\partial M_i|\leq c m_{BY}(\partial M_i,g),$$
where $c$ is a positive constant. Moreover, equality holds for some compact $\partial M_i$ if and only if $M^n$ is isometric to the standard hemisphere $\mathbb{S}^n_+.$
\end{theorem}

We observe that the isometrical embedding assumption in Theorem \ref{massc} was necessary to apply the Riemannian positive mass theorem. As an immmediate consequence is the following result.

\begin{corollary}\label{cmassc}
Let $(M^n,g, f,h)$, $n\geq 3$ be a compact Einstein-type manifold with (possibly disconnected) boundary and positive scalar curvature satisfying $h\leq R_gf/(n-1).$ Suppose that each boundary component $(\partial M_i,g)$ can be isometrically embedded in $\mathbb{R}^n$ as a convex hypersurface. Then, we have
\begin{eqnarray}\label{cormas}
|\partial M_i|\leq \frac{n-1}{n-2}c^2\displaystyle\int_{\partial M_i}(R^{\partial M_i}+|\mathring{II_i}|^2)dS_g,
\end{eqnarray}
for some positive $c$, where $\mathring{II_i}$ is the traceless second fundamental form of $\partial M_i$ as hypersurface of $\mathbb{R}^n$. The equality holds for some connected component of the boundary if and only if $(M^n,g)$ is isometric to the round hemisphere $\mathbb{S}^n_+$.
\end{corollary}

This article is organized as follows. In Section \ref{prel}, we review some basic facts and useful Lemmas on Einstein-type manifold that will be used in the proofs of the main results. Finally, Section 3 collects the proofs of Theorem  \ref{estlam}, \ref{boundlambda1} and \ref{massc} and Corollary \ref{cmassc}.\\

\textbf{Acknowledgement} The author would like to thank Ernani Ribeiro Jr. for his valorous comments about this work.

\section{Preliminaries}\label{prel}
In this section we consider $(M^n, g,f,h)$ be an Einstein-type manifolds as defined in the Section \ref{intro}. Here, we present some basic facts about this manifold that will be useful for the conclusion of our results announced before. Einstein-type manifold with boundary constitutes as a natural extension of static spaces. The following lemma corresponds to Lemma 1 of \cite{andrade2024some} and establishes the connection between these classes of manifolds.

\begin{lemma}[\cite{andrade2024some}]\label{lemaRconst} Let $(M^n, g,f,h)$ be an Einstein-type manifold. The scalar curvature of $(M^n,g)$ is constant if and only if $Rf-(n-1)h$ is constant.
\end{lemma}

Observe that this result suggests that the Einstein-type manifold alone does not implies the constancy of the scalar curvature, which occurs in the case of vacuum static and $V$-static spaces. Moreover, in the proof of this lemma was showed that $\frac{1}{2}f\nabla R=\nabla(Rf-(n-1)h),$ then  since $f=0$ on $\partial M,$ we infer that $Rf-(n-1)h=C,$ where $C$ is a constant, and in particular on $\partial M$ we obtain $-(n-1)h=C$. By other hand, by \eqref{eq1}, we conclude that $h=0$
 on $\partial M.$ This shows that $\partial M$ of an Einstein-type manifold is totally geodesic. In fact, on Einstein-type manifolds with boundary we can suppose that $|\nabla f|$ is a constant non-null on each connected component of $\partial M$. Thus, we can consider the normal vector on $\partial M$ defined by $\nu=-\frac{\nabla f}{|\nabla f|}. $ From now on, consider an orthonormal frame $\{e_i\}_{i=1}^n$ with $e_n=\nu.$ Thus, from the second fundamental form formula, for $1\leq a,b,c,d\leq n-1,$ one obtains that 
$$\beta_{ab}=\langle \nabla_{e_a}\nu,e_b\rangle=-\frac{1}{|\nabla f|}\langle \nabla_{e_a}\nabla f,e_b\rangle=-\frac{1}{|\nabla f|}\nabla_a\nabla_bf=0,$$
and hence, $\partial M$ is totally geodesic. Thus, from Gauss equation, i.e.,
$$R^{\partial M}_{abcd}=R_{abcd}-\beta_{ad}\beta_{bc}-\beta_{ac}\beta_{bd},$$ we deduce that 
$$R^{\partial M}_{abcd}=R_{abcd}.$$ Then, taking the trace in the coordinates $b, d$, we obtain $$R^{\partial M}_{ac}=R_{ac}-R_{ancn}$$ and finally we conclude that the scalar curvature on $\partial M$ is given by 
\begin{equation}\label{Rcurv}
R^{\partial M}=R-2R_{nn}.
\end{equation}

 In  \cite{coutinho2019static} was proved a following divergence formula that  will be a fundamental result in the proof of Theorem \ref{massc}. We observe that if $(M^n, g,f,h)$ is an Einstein-type manifold, it is a simple calculation to verify that$f\mathring{Ric}=\mathring{\nabla^2f}.$ 

\begin{lemma}[\cite{coutinho2019static}]\label{divf}Let $(M^n, g)$ be a Riemannian manifold and $f$ is a smooth function satisfying $f\mathring{Ric}=\mathring{\nabla^2f}$. Then, in the interior of $M,$ one has
$$div\left[\frac{1}{f}(\nabla_g|\nabla_g f|^2-\frac{2\Delta_gf}{n}\nabla_gf)\right]=2f|\mathring{Ric_g}|+\frac{n-2}{n}\langle \nabla_gR_g,\nabla_gf\rangle .$$
\end{lemma}

%\begin{itemize}
%\item Suppose that $|\nabla f|\neq 0$ in each $\partial M.$
%\item $\beta_{ab}=\dfrac{h}{|\nabla f|}g_{ab}\rightarrow \dfrac{H}{n-1}=\dfrac{h}{|\nabla f|}.$
%\item $R-2Ric(\nu,\nu)=R^{\partial M}-\frac{(n-2)}{(n-1)}H^2, \nu= -\nabla f\|\nabla f|.$
%\item $(n-1)(n-2)w_{n-1}^{\dfrac{2}{n-1}}\geq R^{\partial M}|\partial M|^{\frac{2}{n-1}}$
%\end{itemize}

\begin{remark}\label{impor}Using the Bonnet-Myers theorem and Bishop-Gromov comparison theorem we note that if $\partial M$ is an Einstein Riemannian manifold $(n-1)$-dimensional with $R^{\partial M}>0,$ then we claim that
\begin{equation}\label{Rvol}
R^{\partial M} \leq (n-1)(n-2)\omega_{n-1}^{\frac{2}{n-1}}|\partial M|^{-\frac{2}{n-1}}.
\end{equation}
In fact, since $\partial M$ is Einstein with $R^{\partial M}>0,$ then there exists a constant $\delta>0$ such that $Ric^{\partial M}=(n-2)\delta g_{\partial M}$, where $\delta=\frac{R^{\partial M}}{(n-1)(n-2)}.$ Applying Bonnet-Myers theorem we obtain that the diameter of $\partial M$ satisfies $diam_{g_{\partial M}}(\partial M)\leq\frac{\pi}{\sqrt{\delta}}.$ At the same time, by Bishop-Gromov comparison theorem we infer
$$|\partial M|\leq |B_{\frac{\pi}{\sqrt{\delta}}}(p)|\leq |\mathbb{S}^{n-1}|_{g_{\delta}}=\delta^{\frac{-(n-1)}{2}}\omega_{n-1},$$
where $p\in\partial M,$ $g_{\delta}=\delta^{-1}g_{can}$ and $\omega_{n-1}$ is the volume of the standard unitary sphere $\mathbb{S}^{n-1}$. Thus, using that  $\delta=\frac{R^{\partial M}}{(n-1)(n-2)},$ we conclude the claim.
\end{remark}

As mentioned in the Introduction about the Jacobi operator $J$ and we has denoted by $\lambda_1$ as the first eigenvalue of the Jacobi operator. With this notation, we obtain the following proposition.

\begin{proposition}\label{lambda1up}
Let $(M^n, g,f, h)$ be a compact type-Einstein manifold with boundary $\partial M.$ Suppose that the scalar curvature of $\partial M$ is constant. Then,
$$\lambda_1\geq \frac{R^{\partial M}-R_{max}}{2},$$
where $R_{max}=max\{R(p);p\in \partial M\}.$ Moreover, if the equality holds, then $R$ is constant on $\partial M$ and the eigenfunction associated to $\lambda_1$ is constant.
\end{proposition}

\begin{proof} Since the Einstein-type manifold has boundary totally geodesic, then using the integrating by parts, we obtain
\begin{eqnarray*}
-\displaystyle\int_{\partial M}\phi J(\phi)dS_g&=&-\displaystyle{\int}_{\partial M}\phi\Delta_{\partial M}\phi dS_g-\displaystyle{\int}_{\partial M}R_{nn}\phi^2dS_g\\
&=&\displaystyle{\int}_{\partial M}|\nabla_{\partial M}\phi|^2dS_g-\displaystyle{\int}_{\partial M}R_{nn}\phi^2dS_g.
\end{eqnarray*}
From \eqref{Rcurv}, we infer 
\begin{eqnarray}\label{eqs1}
-\displaystyle\int_{\partial M}\phi J(\phi)dS_g\geq \displaystyle\int_{\partial M}|\nabla_{\partial M}\phi|^2dS_g+\left(\frac{R^{\partial M}-R_{max}}{2}\right)\displaystyle{\int}_{\partial M}\phi^2dS_g,
\end{eqnarray}
for all $\phi\in C^{\infty}(\partial M).$ In particular, choosing $\phi$ such that $J(\phi)=-\beta_1\phi$ one sees that
\begin{equation}\label{eqs2}
\lambda_1\geq\frac{\displaystyle\int_{\partial M}|\nabla_{\partial M}\phi|^2dS_g}{\displaystyle{\int}_{\partial M}\phi^2dS_g}+\left(\frac{R^{\partial M}-R_{max}}{2}\right)\geq \frac{R^{\partial M}-R_{max}}{2}.
\end{equation}
If the equality holds in \eqref{eqs2}, then follows directly from \eqref{eqs1} and \eqref{eqs2} the result. Thus, the proof is finished.
\end{proof}

Before to present our estimative for the area of the boundary in terms of the Brown-York mass, we prove a lemma that will be fundamental in the proof of our result. In this way, we shall adapt some arguments in \cite{costa2023geometry, yuan2023brown}. We consider the conformal metric $\bar{g}=v^{\frac{4}{n-2}},$ with
\begin{eqnarray}\label{valpha}
v=(1+\alpha f)^{\frac{-(n-2)}{2}}\ \text{and}\ \alpha^{-1}=\displaystyle{\text{max}_M}\left(f^2+\frac{n(n-1)}{R_g}|\nabla f|^2\right)^{1/2}.   
\end{eqnarray}
We will prove, under certain condition, a necessarily and sufficient for $\alpha$ to be a constant and $\Delta_gf=-\frac{R_g}{(n-1)}f.$
We aid of these observations, we have the following lemma.

\begin{lemma}\label{lemmaconf}
Let $(M^n,g,f, h)$, $n\geq 3,$ be a compact Einstein-type with $R_g>0$ and $h\leq R_gf/(n-1).$ Then, the scalar curvature $R_{\bar{g}}$ with respect to the conformal metric $\bar{g}$ is non-negative, i.e., $R_{\bar{g}}\geq 0.$ Moreover, $R_{\bar{g}}=0$ if and only if  $h= R_gf/(n-1)$ and $f^2+\frac{n(n-1)}{R_g}|\nabla f|^2$ is constant on $(M^n,g).$ In this case, $R_g$ is constant.  
\end{lemma}
\begin{proof}
First, a direct computation, by using the conformal metric define above, we have
\begin{eqnarray*}
\Delta_gv&=&\Delta_g\left((1+\alpha f)^{-\frac{(n-2)}{2}}\right)\\
          &=&\text{div}_g\nabla_g\left((1+\alpha f)^{-\frac{(n-2)}{2}}\right)\\
           &=&\text{div}_g\left(-\frac{(n-2)}{2}(1+\alpha f)^{-\frac{n}{2}}\alpha\nabla_gf\right)\\
           &=& -\frac{(n-2)}{2}\alpha(1+\alpha f)^{-\frac{n}{2}}\Delta_gf+\frac{n(n-2)}{4}\alpha^2(1+\alpha f)^{-\frac{n+2}{2}}|\nabla_gf|^2.
\end{eqnarray*}
Since $h\leq \frac{Rf}{n-1}$ and from \eqref{eq1}, we obtain $\Delta_gf\geq -\frac{fR_g}{n-1}.$ This implies that
\begin{eqnarray}\label{eqnablav}
\Delta_gv\leq \frac{n(n-2)}{4}\alpha(1+\alpha f)^{-\frac{(n+2)}{2}}\left(\frac{2R_gf}{n(n-1)}(1+\alpha f)+\alpha|\nabla_g f|^2\right).
\end{eqnarray}
In another direction, by formulae for conformal metric, the scalar curvature in the metric $\bar{g}$ is given by
$$R_{\bar{g}}=v^{-\left(\frac{n+2}{n-2}\right)}\left(R_g v-\frac{4(n-1)}{n-2}\Delta_gv\right).$$
Then, we deduces that
\begin{eqnarray*}
\Delta_g v&=&\frac{(n-2)}{4(n-1)}\left(R_gv-v^{\frac{n+2}{n-2}}R_{\bar{g}}\right)\\
&=&\frac{(n-2)}{4(n-1)}\left(R_g(1+\alpha f)^{-\frac{(n-2)}{2}}-(1+\alpha f)^{-\frac{(n+2)}{2}}R_{\bar{g}}\right)\\
&=& \frac{(n-2)}{4(n-1)}(1+\alpha f)^{-\frac{(n+2)}{2}}(R_g(1+\alpha f)^2-R_{\bar{g}}).
\end{eqnarray*}
This combined with \eqref{eqnablav} gives
$$R_g(1+\alpha f)^2-R_{\bar{g}}\leq n(n-1)\alpha\left(\frac{2R_g}{n(n-1)}f(1+\alpha f)+\alpha |\nabla f|^2\right).$$
Thus,
\begin{eqnarray}\label{eqRbar}
R_{\bar{g}}&\geq& R_g(1+\alpha f)^2-n(n-1)\alpha\left(\frac{2R_g}{n(n-1)}f(1+\alpha f)+\alpha |\nabla f|^2\right)\\
&=&R_g\left(1-\alpha^2\left(f^2+\frac{n(n-1)}{R_g}|\nabla_gf|^2\right)\right).
\end{eqnarray}
Now, using the value chosen for $\alpha$ in \eqref{eqRbar}, we conclude that $R_{\bar{g}}\geq 0,$ as asserted.

Finally, if the equality holds in \eqref{eqRbar}, i.e., $R_{\bar{g}}=0,$ if and only if $h=\frac{R_gf}{n-1}$ and $\alpha^{-1}=\left(f^2+\frac{n(n-1)}{R_g}|\nabla_g|^2\right)^{1/2}$ is constant over $M.$ In particular, $R_g$ is constant.
\end{proof}

\section{Proof of main results}

In this Section, we shall present the proofs of Theorems \ref{estlam}, \ref{boundlambda1} and \ref{massc} and Corollary \ref{cmassc}.

\subsection{Proof of Theorem \ref{estlam}}
\begin{proof}
By definition of first eigenvalue of the Jacobi operator and using that the boundary of Einstein-type manifold is totally geodesic, we deduce
$$\lambda_1\displaystyle\int_{\partial M}\varphi^2dS\leq -\displaystyle\int_{\partial M}\varphi J(\varphi) dS =\displaystyle\int_{\partial M}|\nabla^{\partial M} \varphi|^2dS -\displaystyle\int_{\partial M}Ric(\nu,\nu)\varphi^2dS,$$
for any $\varphi\in C^{\infty}(\partial M)$. We take $\varphi=1,$ then by \eqref{Rcurv}, we obtain
\begin{eqnarray*}
\lambda_1|\partial M|&\leq& -\frac{1}{2}\displaystyle\int_{\partial M}RdS+\frac{1}{2}\displaystyle\int_{\partial}R^{\partial M}dS\\
&\leq-&\frac{1}{2}R_{min}|\partial M| +\frac{1}{2}\displaystyle\int_{\partial M}R^{\partial M}dS\\
&=-&\frac{1}{2}R_{min}|\partial M| +\frac{1}{2}\displaystyle\int_{\partial M}R^{\partial M}dS\\
&\leq& \frac{1}{2}(n-1)(n-2)\omega_{n-1}^{\frac{2}{n-1}}|\partial M|^{\frac{(n-3)}{n-1}}-\frac{R_{min}}{2}|\partial M|\\
&=&\frac{1}{2}\left(Y(\mathbb{S}^{n-1}, g_{can})|\partial M|^{\frac{(n-3)}{n-1}}-R_{min}|\partial M|\right),
\end{eqnarray*}
where in the fourth line we have used the estimative \eqref{Rvol}. Thus,
$$\lambda_1\leq \frac{1}{2}(Y(\mathbb{S}^{n-1}, g_{can})|\partial M|^{\frac{-2}{n-1}}-R_{min}).$$
Now,  observe that we can rewrite \eqref{Rvol} as $R^{\partial M}\leq Y(\mathbb{S}^{n-1}, g_{can})|\partial M|^{\frac{-2}{n-1}}.$ Thus, if the equality holds in \eqref{lambda_1}, then from the equality case for the Bishop-Gromov theorem, we conclude that the boundary is isometric to a round sphere $\mathbb{S}^{n-1}.$
%$$\displaystyle\int_{\partial M}R^{\partial M}=Y(\mathbb{S}^{n-1}, g_{can})|\partial M|^{\frac{(n-3)}{n-1}}=(n-1)(n-2)\omega_{n-1}^{\frac{2}{n-1}}|\partial M|^{\frac{(n-3)}{n-1}}.$$
\end{proof}

\subsection{Proof of Theorem \ref{boundlambda1}}

We adapt some arguments in \cite{barros2019rigidity} and \cite{diogenes2024geometric}. Since the boundary of Einstein-type manifold is totally geodesic, we obtain that
$$\lambda_1\displaystyle\int_{\partial M}\varphi^2dS\leq -\displaystyle\int_{\partial M}\varphi J(\varphi) dS =\displaystyle\int_{\partial M}|\nabla^{\partial M} \varphi|^2dS -\displaystyle\int_{\partial M}Ric(\nu,\nu)\varphi^2dS,$$
for any $\varphi\in C^{\infty}(\partial M)$. We take $\varphi=1,$ then from \eqref{Rcurv}, we deduce
\begin{eqnarray}\label{eqvoll}
\lambda_1|\partial M|&\leq& -\frac{1}{2}\displaystyle\int_{\partial M}RdS+\frac{1}{2}\displaystyle\int_{\partial}R^{\partial M}dS\nonumber\\
&\leq-&\frac{1}{2}R_{min}|\partial M| +\frac{1}{2}\displaystyle\int_{\partial M}R^{\partial M}dS.
\end{eqnarray}
By assumption, $R^{\partial M}>0,$ then from Remark \ref{impor} combined \eqref{eqvoll} guive us
\begin{eqnarray}\label{lambdaR}
2\lambda_1+R_{min}\leq(n-1)(n-2)\omega_{n-1}^{\frac{2}{n-1}}|\partial M|^{\frac{-2}{n-1}}.\end{eqnarray}
Now, we will prove that $2\lambda_1+R_{min}>0.$ Suppose, by contradiction, that $\lambda_1+R_{min}\leq 0.$ Then, by Proposition \ref{lambda1up}, we infer
$R^{\partial M}-R_{max}\leq 2\lambda_1\leq -R_{min}.$ This shows that $R^{\partial M}< R_{max}$. But, it is a contradiction with our hypothesis. Thus, the stated inequality follows. 

Moreover, if the equality holds in \eqref{lambdaR}, then \eqref{Rvol} also becomes an equality. Therefore, it follows by Bishop-Gromov´s Theorem from the equality case that the boundary is isometric to round sphere $\mathbb{S}^{n-1}.$

%\section{Boundary area estimative in terms of the Brown-York mass}

%In our next result, we will obtain a boundary estimate for compact type-Einstein manifold with (possibly disconneted) boundary in terms of Riemannian Brown-York mass $m_{By}.$ More precisely, we have the following result.

\subsection{Proof of the Theorem \ref{massc}}

\begin{proof} 
First, we shall adapt some arguments in \cite{costa2023geometry,yuan2023brown}. From \eqref{valpha}, we have that
\begin{eqnarray}\label{valpha1}
v=(1+\alpha f)^{\frac{-(n-2)}{2}}\ \text{and}\ \alpha^{-1}=\displaystyle{\text{max}_M}\left(f^2+\frac{n(n-1)}{R_g}|\nabla f|^2\right)^{1/2}.   
\end{eqnarray}
We claim that the mean curvature $H^i_{\bar{g}}$ of $\partial M_i$ with respect to the conformal $\bar{g}=v^{\frac{4}{n-2}}$ is strictly positive. In fact, from \eqref{valpha1} we infer that $v|_{\partial M}=1$, because $f|_{\partial M}=0.$ Moreover, $\bar{g}=g$ over the boundary and $(\partial M_i,\bar{g})$ is isometric to $(\partial M_i, g)$, by hypothesis can be isometrically embedded in $\mathbb{R}^n$ as a convex hypersurface with mean curvature $H_0^i,$ induced by the Euclidean metric. We recall that the mean curvature of $\partial M_i$ with respect to $\bar{g}$ is given by
$$H^i_{\bar{g}}=v^{\frac{-2}{n-2}}\left(H^i_g+\frac{2(n-1)}{n-2}\partial_{\nu}log(v)\right)=\frac{2(n-1)}{n-2}\partial_{\nu}log(v),$$
here $\nu=-\frac{\nabla f}{|\nabla f|}$ and $H_g^i=0.$ It is simple calculation to verify that $\partial_{\nu}log(v)=\left(\frac{n-2}{2}\right)\alpha|\nabla f|$ on $\partial M_i.$ Thus,
\begin{eqnarray}\label{Hgi}
H^i_{\bar{g}}=(n-1)\alpha |\nabla f|.
\end{eqnarray}
This proves that $H^i_{\bar{g}}>0,$ as claimed.

Proceeding, by \eqref{Hgi}, we obtain
\begin{eqnarray}\label{massabar}
m_{BY}(\partial M_i,\bar{g})&=&\displaystyle\int_{\partial M_i}(H_0^i-H^i_{\bar{g}})dS_g\nonumber\\
&=&\displaystyle\int_{\partial M_i}(H_0^i-H^i_{{g}})dS_g+\displaystyle\int_{\partial M_i}(H_g^i-H^i_{\bar{g}})dS_g\nonumber\\
&=&m_{BY}(\partial M_i,g)-\alpha(n-1)|\nabla f||_{\partial M_i}|\partial M_i|.
\end{eqnarray}

 Therefore, by Lemma \ref{lemmaconf}, we infer that $R_{\bar g}\geq 0.$ Thus, we can use the Riemannian Positive Mass Theorem for the Brown-York mass ( see \cite{shi2002positive} and \cite{yuan2023brown}) to conclude that $m_{BY}(\partial M_i,\bar{g})\geq 0.$ This proves that
\begin{eqnarray}\label{massi}
|\partial M_i|\leq \dfrac{1}{(n-1)\alpha|\nabla f|_{\partial M_i}}m_{BY}(\partial M_i,g)
\end{eqnarray}

Next, if equality holds in \eqref{massi} for some component $\partial M_{i_0}$, then necessarily we have 
$$m_{BY}(\partial M_{i_0},\bar{g})=0.$$

Thus, from the equality case in the Positive Mass Theorem for the Brown-York mass, we deduce that the conformal metric
$\bar{g}$ and this implies that $(M^n,\bar{g})$ is isometric to a bounded domain in $\mathbb{R}^n.$ Since $R_{\bar{g}}=0,$ then by Lemma \ref{lemmaconf}, we conclude that $\Delta_gf=-\frac{R_g}{n-1}f$ and $\frac{n(n-1)}{R_g}|\nabla_gf|^2+f^2$ is constant. Consequently, we obtain that $h=fR_g/(n-1)$. Finally, using Lemma \ref{lemaRconst}, we conclude that $R_g$ is constant over $M$. Now, using that $\frac{n(n-1)}{R_g}|\nabla_gf|^2+f^2$ is constant, we infer $\nabla_g|\nabla_g f|^2-\frac{2\Delta f}{n}\nabla_gf=0.$

Now, using the Lemma \ref{divf}, we obtain that $\mathring{Ric_g}=0,$ i.e. $(M^n, g)$ is Einstein with $R_g>0$. Moreover, since  $\partial M$ is totally geodesic, by Proposition 1 in \cite{coutinho2019static}, we conclude that $(M^n, g)$ is isometric to $\mathbb{S}^n_+$.

Reciprocally, if $(M^n, g)$ is isometric to $\mathbb{S}^n_+,$ then $R_g=n(n-1)$ and $\partial M=\mathbb{S}^{n-1}.$ Thus, by involking the Lemma \ref{lemmaconf}, one deduces that  $h=Rf/(n-1)$, i.e., $
\Delta_gf=-nf.$ We will prove that $f^2+|\nabla_g f|^2$ is constant on $M.$ In fact,
\begin{eqnarray*}
\nabla_g(f^2+|\nabla_g f|^2)&=&2f\nabla_gf+2\nabla_g^2f(\nabla_gf)\\
&=&-\frac{2}{n}\Delta_gf\nabla_gf+2\nabla_g^2f(\nabla_gf)\\
&=&(\mathring{\nabla_g^2f})(\nabla_gf)=2f\mathring{Ric_g}(\nabla_gf)=0.
\end{eqnarray*}
Thus, $f^2+|\nabla_g f|^2$ is constant on $M.$ This implies that $\alpha^{-2}=\left.(f^2+|\nabla_g f|^2)\right|_{\partial M}=|\nabla_g f|_{\partial M}^2$, i.e., $\alpha=|\nabla_g f|_{\partial M}^{-1}.$
On the other hand, since $(M^n, g)$ is isometric  $\mathbb{S}^n_+,$ we infer that $\partial M=\mathbb{S}^{n-1}.$ Then, the Brown-York mass of $\mathbb{S}^{n-1}$ is given by
$$m_{BY}(\mathbb{S}^{n-1})=\displaystyle\int_{\mathbb{S}^{n-1}}(n-1)d_{\mathbb{S}^{n-1}}=(n-1)\omega_{n-1},$$
where $\omega_{n-1}$ is the volume of the round $(n-1)$-dimensional sphere. Therefore,
$$\dfrac{m_{BY(\partial M,g)}}{(n-1)\alpha\nabla_gf}=\omega_{n-1}=|\partial M|,$$
which proves the equality in \eqref{massi}.
\end{proof}
\subsection{Proof of Corollary \ref{cmassc}}

\begin{proof}
First, we observe that $R_{\bar{g}}^{\partial M_i}=R_{{g}}^{\partial M_i}$, because we suppose that $(\partial M_i,\bar{g})$ and $(\partial M_i,g)$ are isometric. Now, since $\partial M_i$ is an embedded hypersurface of $\mathbb{R}^n,$ then by Gauss'equation for $\partial M_i,$ we obtain
$$R_{\bar{g}}^{\partial M_i}=(H_0^i)^2-|II_i|^2=\frac{n-2}{n-1}(H_0^i)^2-|\mathring{II}_i|^2,$$
where $II_i$ and $\mathring{II}_i$ denotes the second fundamental form  and traceless second fundamental form of $\partial M_i$, respectively. Using \eqref{massi}, we obtain
$$|\partial M_i|\leq \frac{1}{(n-1)\alpha|\nabla_g f|_{\partial M_i}}m_{BY}(\partial M_i, g)= \frac{1}{(n-1)\alpha|\nabla_g f|_{\partial M_i}}\displaystyle\int_{\partial M_i}H_0^idS_g,$$
where we have used that $H^i_g=0$ and $c= \frac{1}{(n-1)\alpha|\nabla_g f|_{\partial M_i}}.$ This implies, by Holder's inequality,
\begin{eqnarray}\label{massic}
|\partial M_i|\leq c^2\displaystyle\int_{\partial M_i}(H_0^i)^2dS_g=c^2\left(\frac{n-1}{n-2}\right)\displaystyle\int_{\partial M_i}(R_g^{\partial M_i}+|\mathring{II}_i|^2)dS_g.
\end{eqnarray}
If holds the equality in \eqref{massic}, then \eqref{massi} also becomes an equality. Therefore, we infer that $(M^n,g)$ is isometric to round hemisphere $\mathbb{S}^n_+.$
Reciprocally, if $M$ is isometric to $\mathbb{S}^n_+$ with standard metric, then $R_g^{\partial M}=(n-2)(n-1)$ and $\mathring{II_i}=0$. Moreover,
$$\left.\left(f^2+\frac{n(n-1)}{R_g}|\nabla_g f|^2\right)\right|_{\partial M}=|\nabla_g f|_{|_{\partial M}}^2$$
and thus, we conclude
$$\dfrac{1}{(n-1)(n-2)\alpha^2|\nabla_g f|_{|_{\partial M_i}}^2}\displaystyle\int_{\partial M_i}(R_g^{\partial M_i}+|\mathring{II}|^2)dS_g=\omega_{n-1}=|\partial M|,$$
this gives the equality in \eqref{massic}. So, the proof is finished.
\end{proof}

%\bibliography{main.bib}

%\bibliographystyle{acm}
\end{document}